\documentclass[12pt,thmsa]{article}
\usepackage{amsfonts}
\usepackage{amssymb}
\usepackage{graphicx}
\usepackage{amssymb,amsmath}
\newtheorem{teo}{Theorem}[section]

\newtheorem{lema}[teo]{Lemma}
\newtheorem{proposition}[teo]{Proposition}

\newtheorem{lemma}[teo]{Lemma}

\newtheorem{obs2}[teo]{Remark}

\newtheorem{tea}{Theorem}[subsection]

\newtheorem{no2}[teo]{Note}

\newtheorem{no3}[tea]{Note}

\newcommand{\Gal}{{\rm Gal}}

\newcommand{\GL}{{\rm GL}}

\newcommand{\F}{{\mathbb{F}}}
\newcommand{\Q}{{\mathbb{Q}}}

\newcommand{\diag}{{\rm diag}}

\def\timehm{\count31=\time \count32=\count31 \divide\count31 by 60
\number\count31 \multiply\count31 by 60 \advance\count32 by
-\count31 :\ifnum\count32<10 0\fi \number\count32}

\newcommand{\qed}{\hfill\rule{2mm}{2mm}}

\def\mod{\mathop{\rm mod}\nolimits}

\def\ideal#1{<\kern-2pt #1\kern-2pt >}

\begin{document}
\title{{\bf  Remarks on Serre's modularity conjecture}
}
\author{Luis Dieulefait
\\
Dept. d'Algebra i Geometria, Universitat de Barcelona;\\
Gran Via de les Corts Catalanes 585;
08007 - Barcelona; Spain.\\
e-mail: ldieulefait@ub.edu\\
 }
\date{\empty}

\maketitle

\vskip -20mm

\begin{abstract} In this article we give a proof of Serre's conjecture for the case of odd level and arbitrary weight.
 Our proof  does not use any modularity lifting theorem in
  characteristic $2$ (moreover, we will not consider at all characteristic $2$ representations at any step of our proof).\\
The key tool in the proof is a very general modularity lifting result of Kisin,
 which is combined with the methods and results of previous articles
  on Serre's conjecture by Khare, Wintenberger, and the author,
   and modularity results of Schoof for abelian varieties of small conductor.
    Assuming GRH, infinitely many cases of even level will also be proved.
\end{abstract}

\section{Introduction}
Let $p>2$ be a prime, and let $\bar{\rho}$ be an odd, irreducible,
 two-dimensional Galois representation of the absolute Galois group of $\Q$ with Serre's weight $k$ and
  level $N$, with values in a finite field $\F$ of characteristic $p$.
   The ``level", or ``conductor", is defined as in [Se] to be the prime-to-$p$
    part of the Artin conductor, see [Se] for the definition of the weight. \\
We will be mainly interested in the case of representations of odd level, although
 some cases of even level will also be considered, but only cases where ramification
  at $2$ is semistable (in the sense of [Ri2]). \\
For such $\bar{\rho}$, and in particular for all representations of odd level, we will
prove Serre's modularity conjecture (assuming GRH in the cases of even level), i.e.,
we will prove that $\bar{\rho}$ is modular (cf. [Se]). As  is well-known,
 for a given prime $p$ (by suitable twisting) it is enough to consider the case of $k \leq p+1$.
  In all steps of the proof, whenever a residual irreducible representation is considered,
   it will be tacitly assumed that such a twist is performed so that the weight satisfies this inequality.\\

In this article, as in previous articles proving special cases of Serre's conjecture,
 modularity is proved by ``propagation", i.e., by applying the principle of ``switching the residual
  characteristic" (originally applied in [Di3] and [KW1] to prove the first cases of Serre's conjecture)
   to reduce the problem to some other case of the conjecture already solved. This ``switching principle"
    follows from a combination of three results:\\
- Existence of crystalline minimal lifts or lifts with prescribed properties ([Di3], [KW1], [Kh], [KW2])\\
- Existence of (strictly and strongly) compatible systems ([Di2]), and\\
- Modularity lifting results \`{a} la Wiles (Taylor-Wiles, Skinner-Wiles, Diamond, Savitt, Kisin).\\
 At this point, the main ``obstacle for propagation" is due to the technical conditions needed for
  the application of these modularity lifting results. Still, let us recall that in several cases,
   like the crystalline of small weight case ($k < p$, assuming that $p \neq 2k-3$, or the
    representation is semistable), or weight $2$ semistable case, it is known that by
     combining different modularity lifting results the lifting is modular without
      imposing any condition on the residual representation, just modularity or
       reducibility (cf. [Di1] for the weight $2$ case and [DM], [Kh] for higher weights).\\

By a suitable combination of ``switchings", using mainly modularity lifting
 results of Kisin (cf. [Ki4]), we will show how the proof of the general
  odd level case can be reduced to the proof of the level $3$ case, a
   case that we will reduce in turn to some modularity results of Schoof
    for semistable abelian varieties of small conductor (cf. [Sc2]). Also, some cases
     of even level will be solved assuming GRH.\\

Since the most important step of the proof is the reduction to the
level $3$ case, let us briefly explain how this is achieved:
  via weight $2$ lifts we ``reduce" the proof to a case where ramification conditions
   are much better so that the result of Kisin in [Ki4] can be applied and by a simple
    process of ``iterated killing ramification" we end up in a case with
     just one prime in the level. The preparatory phase that precedes the killing
      ramification process is thus a key step. It is based on the observation that
       one can {\bf transfer ramification} from the given set of primes in the level to
        other sets of sufficiently large primes (using weight $2$ compatible systems), so by choosing {\bf special}
         characteristics, for example primes that are pseudo Sophie Germain,
          one can gain {\bf full  control} on the type of ramification. \\
          This way, after the preparatory phase, we are able to
          check, relying heavily on results of Caruso (cf. [Car]),
          that the required conditions for the application of
          Kisin's modularity lifting theorems are satisfied and thus
          we can
          perform the sought for ``killing ramification": this way
          of attacking the conjecture was proposed, in naive form,
          by Khare-Wintenberger in their first joint paper [KW1], nevertheless they follow a different path
           in their solution of
           Serre's conjecture (see [KW2] and [KW3]) that is based purely on ``weight reduction"
            and depends strongly on proving new modularity
           lifting theorems in characteristic $2$. \\
           Our results imply in particular that for applications to
           residual modularity, the main theorem in [Ki4], combined
            with
           previous modularity lifting theorems, is universal in
           the sense that we can always ``reduce" our problem to a
           situation where the conditions required to apply this
           theorem will always be satisfied.\\
           In conclusion, we believe that the proof that
           ``killing ramification" is possible is key to the theory, it is clearly of independent interest and
           in particular it is a result that will very likely
           have applications in generalizations of the conjecture. \\
           
           In fact, the author has recently (cf. [Di4]) obtained a proof of non-solvable base change for $\GL_2$, and the proof is very much related to the strategy used in this paper to ``propagate modularity", mainly via the main theorem in [Ki4]. In particular the results of Caruso are also needed at some steps of the proof. The method consists on linking, via a suitable chain of congruences, any modular Galois representation with some other modular Galois representation of small level and weight such that it can be base-changed, and at any step of the chain propagation of modularity will be ensured after checking the conditions needed to apply some modularity lifting theorem, mainly the one in [Ki4]. It is by no means clear that a proof of this result could be obtained following the longer row of ``weight reductions" instead of the ``killing ramification" strategy similar to the one in this paper.\\



Remark: One of the special features of the proof presented in this
paper is that it does not require modularity lifting theorems in
characteristic $2$.
 This can be particularly helpful when dealing with generalizations of Serre's
 conjecture to totally real number fields where $2$ is ramified, such as $\Q(\sqrt{2})$ (keeping the assumption that the conductor of the representation is odd).
  Due to the fact that, together with the usual technical difficulties to work in characteristic $2$,
   several modularity lifting results extend to totally real number fields but only with the
    assumption that the prime $p$ is unramified, the possibility to avoid working in
     characteristic $2$ may be key when proving cases of Serre's conjecture over these fields. In fact, as in [Di4], it is sometimes key to be able to avoid all characteristics that do ramify in a given number field, and this is something that clearly can not be done using the method of ``weight reduction".\\
     Also, it will be transparent to the reader that the present
     proof, if combined with modularity lifting theorems in
     characteristic $2$ as those used in the proof of Serre's
     conjecture obtained by Khare and Wintenberger, would extend to
     the case of arbitrary even conductor, but it will keep its
     specific features: it will be a more direct proof by ``killing
     ramification" instead of a proof purely by weight reduction.\\

     The results in this paper were obtained in the period March-May
     2006.\\

 Acknowledgments: I want to thank Jorge Jimenez, Jacques Tilouine, Jean-Pierre Serre and Xavier Caruso
 for useful comments. I want to thank Ren\'{e} Schoof for explaining me his results.
  Also, I want to mention that my attendance to the conference ``p-adic representations, modularity, and beyond" at the American
   Institute of Mathematics (Palo Alto, California, USA) on February 2006, was key for this work. Therefore, I want
     to thank the staff of AIM, and specially the organizers of the conference,
     K. Kedlaya and D. Savitt.\\
  Finally, I would also like to thank the referee for his very useful and detailed comments.\\

  We would like to conclude this introduction by recalling to the reader that this paper and our previous work on Serre's conjecture, as well as the work of Khare and Wintenberger, rely heavily on the potential modularity results of R. Taylor (cf. [Ta1], [Ta2]): the role of potential modularity is fundamental at several key steps of the proof (as can be seen in [Di2], [Di3],  [KW1], [Kh] and [KW3]) and it has also been (at least to me) inspirational.

\section{Preliminary results}

$\bullet$ Crystalline minimal lifts and strongly compatible systems:\\

Let $\bar{\rho}$ be a residual representation of weight $k$, odd level $N$, and odd characteristic $p$. Assume that $2 \leq k \leq p+1$, and the image of $\bar{\rho}$ is non-solvable.\\

For such a representation $\bar{\rho}$, we will need to consider a ``crystalline minimal lift" $\rho_p$ as defined and constructed in [Di3] and in [KW1] (except for the case $k=p$ solved in [KW3]), and the strongly compatible system (in the sense of [Ki2]) $\{ \rho_\ell \}$ containing it (whose existence was proved in [Di2]). The reader should be aware that this kind of lift is usually called ``minimal lift" in the references.\\

We recall the definition of such a lift and system: $\{ \rho_\ell \}$ is a system of de Rham Galois representations of Hodge-Tate weights $\{0, k-1 \}$, if $\ell \nmid N$ the representation $\rho_\ell$ is crystalline (note that this is the case, in particular, for $\ell = p$), for every odd $\ell$ the ramification set of $\rho_\ell$ is the union of the set of primes dividing $N$ and $\ell$, the system is strongly compatible and the representation $\rho_p$ is a lift of $\bar{\rho}$ that is minimally ramified at all primes other than $p$.\\

$\bullet$ Weight $2$ minimal lifts and strongly compatible systems:\\

Let $\bar{\rho}$ be a residual representation of weight $k$, odd level $N$, and odd characteristic $p$. Assume that $ 2 < k \leq p$, and the image of $\bar{\rho}$ is non-solvable.\\

At several steps of the proof we will need to consider for such a $\bar{\rho}$ a p-adic ``weight $2$ minimal lift" $\rho_p$ defined as in [KW2], theorem 5.1(2), and to introduce $\rho_p$ in a strongly compatible system  $\{ \rho_\ell \}$. The proof of the existence of such a lift and compatible system is given in [KW3].\\

The description of the representations in this compatible system is the following: First of all, the lift $\rho_p$ is a p-adic lift of $\bar{\rho}$ that is minimally ramified at every prime different from $p$. Then, the system is strongly compatible and every representation  $   \rho_\ell \in \{ \rho_\ell \}$ with $\ell$ odd ramifies only at primes dividing $N$, $p$ and $\ell$, it is Barsotti-Tate if $\ell \nmid p N$ and potentially semistable of weight $2$ or potentially Barsotti-Tate for any odd $\ell$ (in other words: it is a system of de Rham representations of Hodge-Tate weights $\{0, 1 \}$), and it has inertial Weil-Deligne parameter at $p$ equal to $(\omega_p^{k-2} \oplus 1, 0)$ where $\omega_p$ is the Teichmuller lift of the $\mod$ $p$ cyclotomic character. In particular
$\rho_p$ is potentially Barsotti-Tate (and Barsotti-Tate over a subextension of $\Q_p(\mu_p)$).\\

$\bullet$ Weight $2$ minimal lifts when $k = p+1$:\\

Another case of weight $2$ minimal lift, with similar assumptions except that we assume that $\bar{\rho}$ has weight $k = p +1$, will also be needed. The existence of such lifts and the strongly compatible system containing it is proved in [KW1] and [Kh]. Here the lift and the system are again ``of weight 2" (i.e., de Rham representations of Hodge-Tate weight $\{0, 1  \}$), the lift is again minimally ramified away from $p$ and the system is strongly compatible, and the difference is that now the inertial Weil-Deligne parameter at $p$ of all representations in the system is of the form $(\rm{id}, N  )$ with $N$ a non-zero nilpotent matrix. In other words, the compatible system has ``semistable" ramification at $p$.\\

 $\bullet$ Raising the level and good-dihedral primes:\\

 Because of the technical conditions needed to apply Kisin's modularity lifting results, in the general case we need to assume that  the residual representations that we encounter through the proof have non-solvable image.\\
 This is handled by a  trick created in [KW2]: after switching, using a compatible system, to a suitable (and sufficiently large) characteristic $t$ they add some extra prime $q$ (sufficiently large) to the level in order that all the representations that we encounter are locally ``good dihedral at $q$" and thus have non-solvable image. The precise references are as follows: see section 2 of [KW2]  for the definition of good-dihedral prime, Theorem 3.4 for the theorem asserting that the proof of Serre's conjecture can be reduced to cases having a good-dihedral prime in the level, Lemma 6.3 for the fact that a good-dihedral prime implies residual big images in all relatively small characteristics, and section 8.4 for the proof of the existence of a suitable characteristic $t$ where we can add the good-dihedral prime $q$ to the level, thus proving Theorem 3.4. \\

The possibility to introduce ``extra ramification" follows from the existence of non-minimal $p$-adic lifts of certain type (non-minimal at $q$), as stated in Theorem 5.1 (4) of [KW2], whose proof is given in [KW3].
As in the case of crystalline minimal lifts (cf. [Di3] and [KW1]) or weight $2$ minimal lifts, the existence of such a lift (which mimics a result that is well-known for modular forms, namely, a case of ``raising the level") follows by combining potential modularity  with the results of Boeckle on presentations of universal deformation rings. \\

 This way, after having ``raised the level", we can assume that at each step we will encounter residual representations having a ``good-dihedral" large prime $q$ in the level, i.e., we have a $\bar{\rho}_p$ and a prime $q$ such that $q^2 \mid N$ with:

 $$ \bar{\rho}_p|_{I_q} = \diag(\psi , \psi^q) \qquad \quad (I) $$
 where the character $\psi$ has order $t^\alpha$, $ t \mid q+1$ an odd prime, with $q$ and $t$ sufficiently large.\\
 This ramification at $q$ will be preserved and the residual image will be non-solvable in all the steps of the proof where the residual characteristic is ``relatively small" if the primes $q$ and $t$ are chosen as in the definition of good dihedral prime, as shown in [KW2], Lemma 6.3.\\
 As a matter of fact, we have to modify slightly the definition of good dihedral prime, because we want to work also in several  characteristics up to a certain bound $B$, a bound certainly larger than the weight and the primes in the level\footnote{during section 4 several auxiliary primes will be introduced, and the bound $B$ can be taken to be any number greater than all these auxiliary primes: in the remark at the end of section 4 we will stress this point again}. Thus, to ensure that also in these characteristics the ramification at $q$ is preserved 
 we modify the definition given in [KW2], Def. 2.1, as follows: the ramification at $q$ is given again by formula (I) but we want to:\\
 i) assume $t$ is greater than all prime factors of $N$ (except $q$), greater than $k$, and greater than $B$.\\
 ii) assume that $q \equiv 1 \pmod{8}$, and $q \equiv 1 \pmod{r}$ for every prime $r$ up to the maximum of: the prime divisors of $N$ different from $q$, $k$, and $B$.  \\

 Remark: this ``modified" definition does not affect at all the proof of existence of such a $q$ and $t$.\\

 As in section 8.4 of [KW2], at a large characteristic $t$ we add the extra ramification at $q$ in order to reduce the proof of Serre's conjecture to representations being locally good-dihedral at $q$. In all steps of the proof we will work in characteristics which are sufficiently small with respect to $q$ and $t$, namely smaller than a certain bound $B$ previously given. Thus, we know that all residual representations that we encounter maintain the good dihedral prime $q$ in their levels and thus have non-solvable image (cf. Lemma 6.3 of loc. cit.). Let us formally summarize these results in the following:
 
 \begin{lema} (Khare-Wintenberger):
\label{teo:good-dihedral} 
Given a Galois representation $\bar{\rho}_p$ in characteristic $p$ ($p$ odd) of Serre's level $N$, weight $k$ and non-solvable image, and a bound $B$ greater than the maximum of $N, p$ and $k$, there exist primes $t$ and $q$ as in the definition of ``locally good-dihedral at $q$ with respect to the bound $B$" given above such that the modularity of $\bar{\rho}_p$ can be deduced from the modularity of a residual representation $\bar{\rho}'_p$ in characteristic $p$ of level $Nq^2$ which is locally good-dihedral at $q$.\\ 
 To establish the modularity of such a $\bar{\rho}'_p$ we are going to proceed through strongly compatible systems which are going to be linked in suitable characteristics, and as long as we work exclusively in characteristics smaller than $B$ all residual representations are going to be locally good-dihedral at $q$ and will thus have non-solvable images.
\end{lema} 

   \section{A very general modularity lifting result of Kisin}

       We will rely heavily on the modularity lifting result of Kisin in [Ki4].  This result will ensure that modularity is preserved in many of the links between compatible systems that we will consider in the proof.
       
       \begin{teo} (Kisin):
       \label{teo:superKisin} Let $\rho$ be a representation of the absolute Galois group of $\Q$ with values on a finite extension of ${\Q}_p$ that is $2$-dimensional, continuous, odd, absolutely irreducible and ramified at finitely many primes. Assume also that  the representation is de Rham of Hodge-Tate weights $\{0, k-1 \}, \; k \geq 2$, and that its residual image is non-solvable. Let us also assume that the following technical condition is satisfied: \\
         $\bar{\rho}|_{D_p}$ is not  isomorphic to an extension of $1$ by the mod $p$ cyclotomic character $\chi$, nor to a twist of such a representation. \\
         Then, if $\bar{\rho}$ is modular, $\rho$ is also modular.
       \end{teo}

We will refer to the technical condition on the residual representation locally at $p$ as ``condition (T)".\\

   We include, for the reader's convenience, the following ``road map" that describes  schematically all the moves that we will do in the proof of Serre's conjecture for odd level.

   \subsection{Road Map of the proof}

   The proof consists of a series of {\bf reductions}: by performing suitable {\bf moves} we reduce the proof in each step to a simpler case.\\

   {\bf Notation}: $\bullet$ We will denote by $\mathcal{S}_k(N)$ the truth of Serre's conjecture in the case of Serre's weight $k$ and level $N$, for any characteristic $p \geq k-1$.\\
  $\bullet$ We indicate the steps of the proof, i.e., the reductions, by a reverse implication:
   $$ \mathcal{P}  \Leftarrow  \mathcal{P}' $$
   meaning that the truth of the case $\cal{P}$ in the left hand side is reduced to the truth of the cases included in proposition $\cal{P}'$ on the right hand side.\\

   With this notation, let us summarize the proof. $N$ denotes an odd integer, $k \geq 2$, and whenever we write a factor $q^2$ in the level this means that $q$ is a good-dihedral prime. In sections 4 to 6 we give the following proof:\\

   By an application of Lemma \ref{teo:good-dihedral}, we have:
   $$ \mathcal{S}_k(N)  \Leftarrow \mathcal{S}_k(N q^2) $$
   In section 4 we prove the following:
   $$ \mathcal{S}_k(N q^2)  \Leftarrow \mathcal{S}_{k'}( q^2), \; \forall k' < B \qquad \quad (3.1)$$
   for a bound $B = B(k,N)$ as in the definition of $q$ being a good-dihedral prime.\\
   In section 5 we prove:
   $$ \mathcal{S}_k(q^2)  \Leftarrow \mathcal{S}_{k'}(3) \; \mbox{and} \;  \mathcal{S}_{k'}(1), \; \forall k'  \qquad \quad (3.2)$$
   Since the level $1$ case was already solved in [Kh], it remains to deal with the level $3$ case.\\
   In section 6 we prove:
   $$ \mathcal{S}_k (3) \Leftarrow \mathcal{S}_{k'}( 3 ), k' =2,4,6$$
   and we show that these small weight cases follow (after switching to characteristic $5$ and taking a minimal weight $2$ lift) from modularity results of R. Schoof for abelian varieties of small conductor.\\

   Remark: Step (3.1) is the most complicated one. In fact, we will first have to transfer twice the ramification to other sets of primes in order to be able to perform the ``killing ramification", i.e., to eliminate all primes in the level except $q$ and reduce to a case of level $q^2$.       Step (3.2) decomposes as follows:
   $$ \mathcal{S}_k(q^2)  \Leftarrow \mathcal{S}_{k'}(q^2), k'= 2, 4  \Leftarrow \mathcal{S}_2(u q^2), u=1,3 \Leftarrow $$
   $$ \Leftarrow \mathcal{S}_2(u q), u= 1, 3  \Leftarrow \mathcal{S}_{k''}(u), u=1,3, \; \forall k'' $$

    \section{Iterated killing ramification: the Sophie Germain trick}

     The first step can roughly be described as follows: Killing ramification at primes in the level, one after the other (using crystalline minimal lifts to go from one to the other) we reduce to a case with just one prime (the good-dihedral prime) in the level.\\

     Since we need to ensure that modularity propagates well in all ``moves" (i.e., every time that we switch the residual characteristic), we have to ensure that the conditions of ``non-solvable residual image" and Condition (T) needed to apply Kisin's result are satisfied (cf. Theorem \ref{teo:superKisin}). Thus, before starting the killing ramification process, we will perform two sets of moves, via weight $2$ lifts, to reduce to a situation where these conditions are satisfied.\\
     We start by applying Lemma \ref{teo:good-dihedral}: by adding extra ramification at a very large prime $q$, given by some character of very large order $t$, as long as we work in characteristics smaller than a certain bound $B$ all representations will be good-dihedral at $q$ and will have non-solvable residual image.\\

     Then, we can assume that we start with a residual representation such that the primes in the level (all odd) are:
     $$ p_1, p_2,......p_r, q$$
     of weight $k$ in characteristic $p$ such that for certain bound $B$ larger than $p, k$ and all the $p_i$ the prime $q$ is a good dihedral prime with respect to $B$ (the bound $B$ should also be larger than  the auxiliary primes $b_i$ and $q_i$ used during the proof).\\
First we transfer the ramification to $r+1$ primes which are ``pseudo Sophie Germain", namely, primes $b_i$ such that $(b_i-1)/2$ is an odd prime or the product of two odd primes. It is proved in [Xi] and [Co] that there are infinitely many such primes, and 
from the proof  it is clear that every prime factor of $(b_i -1)/2$ is going to grow with $b_i$. In fact, it is a common feature of all proofs using sieve theory that primes or semiprimes are obtained by restricting to numbers which are relatively prime to all primes up to some increasing function $f(x)$ of $x$, so when $x$ goes to infinity the result of Coleman and Xie implies not only that the number of pseudo Sophie Germain primes  up to $x$ goes to infinity but also that the odd primes dividing the pseudo Sophie Germain primes minus $1$ also grow with $x$, thus go to infinity\footnote{we thank Jorge Jimenez for these clarifications on sieve theory. He also
 reconstructed in full detail the proof of Xie and Coleman for the specific case of pseudo Sophie Germain primes, and in his proof 
  the functions $x^{1/8}$ and $x^{1/3}$ are used, showing in particular that both odd prime factors are at least as large as
   $x^{1/8}$}. \\
 So we assume that both the primes $b_i$ and the odd prime factors of $b_i - 1$ are larger than all the $p_i$ and than $p$ and $k$.\\
Assume also that $b_1 < b_2 < ..... < b_{r+1}$. \\

\begin{lema}
\label{teo:Sophie} If a residual Galois representation $\bar{\rho}_p$ has Serre's level containing only the primes $p_1, ..., p_r$ and $q$, Serre's weight $k \leq p+1$, non-solvable image, and  is locally good-dihedral at $q$ ($q$ larger than some bound $B$), given $r+1$ pseudo Sophie Germain primes $b_1 < b_2 < ..... < b_{r+1}$ as above (larger than the primes $p_j$, than $p$ and $k$ but smaller than $B$), the modularity of $\bar{\rho}_p$ follows from the modularity of a strongly compatible system $\{ \rho'_\ell \}$ of de Rham representations of weight $2$ (i.e., of Hodge-Tate weights $\{0 , 1  \}$) whose ramification set is:
$b_1, b_2, ..... , b_{r+1}$  and $q$, and such that the inertial Weil-Deligne parameter at each $b_i$ is either of the form $(id, N)$ or $( \omega_{b_i}^{(b_i -1)/2} \oplus 1, 0)$. Furthermore, for every prime $p' < B$, the $p'$-adic representation $ \rho'_{p'}$ in this compatible system will be residually locally good dihedral at $q$.
\end{lema}

\bf{Proof}: \rm We proceed to transfer ramification to these $r+1$ primes, via weight $2$ lifts: starting in characteristic $p$, take a crystalline minimal lift and move to characteristic $b_1$, and reduce mod $b_1$. Take a minimal weight $2$ lift
and move to characteristic $p_1$, and reduce mod $p_1$. Here, as usual (and the same in all future steps), we twist this residual representation if necessary in order to assume that its Serre's weight is at most $p_1 + 1$. Then take a crystalline minimal lift and  move to characteristic $b_2$, reduce mod $b_2$, take a minimal weight $2$ lift, move to $p_2$, and so on, at the end we take a minimal weight $2$ lift in characteristic $b_{r+1}$.\\
     Modularity is preserved in all these moves due to Modularity Lifiting Theorems (M.L.T.) of Kisin and Diamond (cf. [Ki1] and [Dia]) in the potentially Barsotti-Tate or weight $2$ semistable case (where here potentially means over any extension), and for crystalline representations of small weight by results of Diamond-Flach-Guo (cf. [DFG]). This means that whenever we are linking a given compatible system $\{ \rho_\ell \}$ with a new compatible system $\{  \rho'_\ell \}$ through a congruence, modularity of the first compatible system can be deduced from the modularity of the second. In the above process, since we are working in characteristics smaller than $B$, it follows from Lemma \ref{teo:good-dihedral} that the good-dihedral prime $q$ ensures that the residual images are non-solvable, and this condition is enough to make the above mentioned M.L.T. work. Notice also for the application of the M.L.T. in [DFG] that when we move to characteristics $b_i$ we are doing so with a strongly compatible system that is unramified at $b_i$, so that its $b_i$-adic member is crystalline, and  the system has non-zero weight $ k_i \leq p_{i-1} +1$ (with $p_0 := p$), which is smaller than $b_i$.\\
     We have transfered the ramification at the $p_i$ to ramification at the $b_i$ introduced in the weight $2$ lifts, given by  characters $\omega_{b_i}^{k_i-2}$, i.e, we have now compatible systems ramified at the primes $b_i$ and at $q$, whose inertial Weil-Deligne parameter at each $b_i$ is $(\omega_{b_i}^{k_i-2} \oplus  1, 0)$.\\
     Observe that this ramification at $b_i$ is introduced when taking a minimal weight $2$ lift in characteristic $b_i$, and it has not changed in the subsequent moves because we are working in characteristics that do not divide $b_i -1$ (by construction, the primes $p_j$ do not divide $b_i-1$, and the sequence $b_i$ is increasing) thus a fortiori they don't divide the order of the character $\omega_{b_i}^{k_i-2}$. We thank the referee for this remark.\\

     For each $b_i$ we move to the odd characteristics (one or two) dividing $b_i -1$, always via weight $2$ systems, so that in these characteristics we kill part of the ramification at $b_i$:  assume first that $(b_i -1)/2 = a_i$ is prime, then in characteristic $a_i$ the nebentypus $\omega_{b_i}^{k_i-2}$ which is a character of order $a_i$ or $2a_i$ (depending on the parity of $k_i$) either becomes a character of order $2$, or the residual representation is unramified or semistable at $b_i$. If $(b_i-1)/2$ is the product of two primes $a_i$ and $a'_i$ then by first moving to characteristic $a_i$, considering the residual representation and taking a minimal lift of it, we change the order of the character giving ramification at $b_i$, now it is a character of order $a'_i$ or $2a'_i$ and so we see that by moving to characteristic $a'_i$ and considering the residual representation  we are again in a case where the ramification at $b_i$, if any, is either semistable or given by a quadratic character.\\
     We conclude that we have linked our given representation to a compatible system such that its $b_i$-adic members, locally at $b_i$, are either semistable or quadratic-crystalline, i.e., crystalline when restricted to a quadratic extension. The primes $b_i$ must be chosen sufficiently distant from each other (in a sense to be specified in the next argument). 
     \qed
     \\
\begin{lemma} Let $\{ \rho'_\ell \}$ be  a strongly compatible system as in the conclusion of the previous lemma. Then, if for every prime $b_i$ in the level where the inertial Weil-Deligne parameter is $(id, N)$  we choose a prime $q_i$ such that $q_i -1$ is divisible by $2 \cdot (b_i - 1)$, all of them smaller than the bound $B$, and we rename the remaining $b_i$ as $q_i$ to ease notation, the modularity of the system can be deduced from the modularity of another strongly compatible system $ \{\rho''_\ell \}$ of weight $2$ having ramification set  $q_1 , q_2, ..... , q_{r+1}$ and $q$ and inertial Weil-Deligne character at $q_i$ equal to $(\omega_{q_i}^{m_i} \oplus 1, 0)$ with exponent $m_i = b_i -1$ for those index $i$ where $q_i \neq b_i$ and $m_i = (b_i-1)/2 = (q_i - 1)/2$ for the others.  Observe that in both cases the exponents grow with $b_i$ and they divide the order of $\omega_{q_i}$. Again, for every prime $p < B$, the $p$-adic representation $ \rho''_{p}$ in this compatible system will be residually locally good dihedral at $q$.
\end{lemma}

\bf{Proof}: \rm
We want to perform a set of moves that allows us to transform all semistable ramification in potentially crystalline ramification. For each $b_i$ such that ramification at it is semistable, we choose a larger prime $q_i$ such that $q_i -1$ is divisible by $2 \cdot (b_i - 1)$. Then we start with a weight $2$ system, we switch to characteristic $b_i$ and reduce mod $b_i$ (if this residual representation has weight $2$ we have eliminated ramification at $b_i$, if not we continue) and take a crystalline minimal lift, which has weight $b_i + 1$, i.e., Hodge-Tate weights $\{0, b_i \}$. We move to characteristic $q_i$, reduce mod $q_i$, and take a weight $2$ lift, thus ramification at $q_i$ is given by the character $\omega_{q_i}^{b_i-1}$. Namely, we have transfered semistable ramification at $b_i$ into potentially crystalline ramification at $q_i$, but most importantly, since the exponent $b_i-1$ divides the order $q_i-1$ of the character, ramification at $q_i$ is crystalline over a {\bf proper} subfield of the cyclotomic field, a subfield $F_i$ such that the degree of the cyclotomic field over it grows with $b_i$ (*). \\
At the end, both at those $b_j$ where ramification was quadratic-crystalline and at  the $q_i$ just considered, we conclude that we have reduced to a case where ramification is always potentially crystalline, over a proper  subfield $F_i$ of the cyclotomic field, as in (*). Just to ease the notation, we rename the primes $b_j$ where ramification was quadratic crystalline as $q_j$. \\
Modularity is preserved through these moves as follows again by an application of the aforementioned M.L.T.
\qed
\\
$$ $$
For the next steps to work, the primes $b_i$ must be assumed to be sufficiently distant from each other so that we can ``interpolate" the primes $q_i$. To be more precise, the primes $b_i$ and $q_i$ should be taken in the following way:
choose first $b_1$ and $q_1 > b_1$, then $b_2 >> q_1$ and $q_2 > b_2$, then $b_3 >> q_2$ and $q_3 > b_3$, and so on (except for $q_j = b_j$ in the quadratic crystalline case). Here we give to the symbol $>>$ the following special meaning: we require the following inequality: $b_i -1  > q_{i-1} + 1$ except for the quadratic crystalline case where  $b_i = q_i$ and we require $(b_i - 1)/2 > q_{i-1} +1 $. \\
It is obvious that we can choose pseudo Sophie Germain primes $b_i$ and primes $q_i$ with $q_i -1$ divisible by $2 \cdot (b_i - 1)$ satisfying these inequalities (just because there are infinitely many pseudo Sophie Germain primes and because of Dirichlet's theorem for primes in arithmetic progressions). \\
 As the referee has remarked, we can also suppose (using Dirichlet's theorem again) that the primes $q_i$ are chosen in such a way that whenever $i \neq j$, $q_i - 1$ is not divisible by $q_j$. This assumption is helpful because it implies that in the following steps we are not going to change the inertial Weil-Deligne parameter at any $q_j$ when working in characteristic $q_i$, for any $i \neq j$. \\
 
 \begin{proposition} Let $ \{\rho''_\ell \}$ be a strongly compatible system as in the conclusion of the previous lemma, whose conductor is divisible only by $q_1 , q_2 , ...., q_{r+1}$ and $q$. Assume also that the primes $b_i$ and $q_i$ have been chosen as in the previous explanation, i.e., not only $q_i -1$ divisible by $2 \cdot (b_i - 1)$ (in cases where $b_i \neq q_i$) but also the  inequalities $b_{i} >> q_{i-1}$ are satisfied and no prime $q_j$ divides $q_i -1$ for any $i, j$. Also, all primes $b_i$ and $q_i$ are smaller than the bound $B$ in the definition of good-dihedral prime (see section 2). Then, we can perform ``killing ramification" in order to link the system to a residual Galois representation of weight and characteristic both smaller than $B$, level $q^2$, and locally good-dihedral at $q$, in such a way that modularity is preserved, in other words, the modularity of the given system follows from the modularity of this residual representation.
 \end{proposition}

\bf{Proof}: \rm At this point, we have to perform iterated killing ramification to eliminate, one after the other, ramification at each of the primes in the level except for the very large good-dihedral prime $q$. Just one warning: this must be done in increasing order. So we begin with a weight $2$ system and we move to $q_1$, the smallest prime in the level, take the residual representation (and, as usual, twist to obtain minimal weight), then a  crystalline minimal lift, move to $q_2$, and so on. After killing ramification at all the $q_i$ we end up with a residual representation of weight smaller than $B$,
of conductor $q^2$ and locally good-dihedral at $q$. Whenever we switch to a prime $q_i$ in the level with a system of weight larger than $2$, we have to check that condition (T) in Kisin's  Theorem \ref{teo:superKisin} is satisfied (in the weight $2$ case the result of [Ki1] applies instead). Since ramification at each $q_i$ is potentially crystalline over a proper subfield $F_i$ of the cyclotomic field of  index $b_i-1$ (as in (*)) or $(b_i -1)/2$ in the quadratic crystalline case (where $b_i = q_i$), and we switch to characteristic $q_i$ with a system of weight at most $q_{i-1} +1$,  we will show in the following lemma that condition (T) is satisfied because we have  $b_i >> q_{i-1}$ (see the paragraph previous to the proposition for the definition of the symbol $>>$), thus concluding the proof of the proposition.  \\

\begin{lema}
\label{teo:caruso} Let $\bar{\rho}_{q_i}$ be a residual representation with non-solvable image as above, i.e., it has a lift corresponding to a potentially crystalline representation of Hodge-Tate weights $\{0,  k \}$ with $k \leq q_{i-1} + 1$ and with inertial Weil-Deligne parameter at $q_i$ equal to  
$(\omega_{q_i}^{m_i} \oplus 1, 0)$ with exponent $m_i = b_i -1$ or $(q_i-1)/2$. In the first case, we suppose that $q_{i-1} + 1 < b_i -1$ and in the second case we suppose that $q_{i-1}  + 1 < (q_i - 1)/2$. Then, the residual representation $\bar{\rho}_{q_i}$ does satisfy the technical condition (T) in Kisin's Theorem \ref{teo:superKisin}.
\end{lema}

\bf{Proof}: \rm We have a potentially crystalline $q_i$-adic representation $\rho$ whose Hodge-Tate weights are $0$ and $k-1$ with $k \leq q_{i-1} + 1$. This representation becomes crystalline over a subfield $F_i$ of the cyclotomic field of ramification degree over $\Q_{q_i}$ equal to $e =
 (q_i -1)/ (b_i -1)$ (which is by construction an even number) or $e=2$ in the second case. Because of the assumed inequalities, in both cases we have:
 $$ (k-1) \cdot e     \leq   q_{i-1}    \cdot e   < q_i -1 $$
 In this situation, we can easily verify that Kisin's condition (T) is satisfied by considering the exponents of tame inertia on the restriction to $F_i$ of the (local at $q_i$) residual representation: if we call these exponents $n_1$ and $n_2$ these exponents are known to satisfy (cf. [Car]):
 $$ 0 \leq n_1 , n_2  \leq e \cdot  (k-1) $$
 If we assume that we are in a case where the action of tame inertia at $q_i$ on $\bar{\rho}$ is given by the characters $\chi^{a}$, $\chi^{a+1}$ for some integer exponent $a$ then we deduce that when considering the restriction of the residual representation to $F_i$ the action of tame inertia is also reducible and the exponents $n_1$ and $n_2$ are both divisible by $e$, so we have $n_1 = e \cdot t_1$ and $n_2 = e \cdot t_2$, and $t_1, t_2$ satisfy:
 $$ 0 \leq t_1 , t_2  \leq  k-1     \qquad \qquad  (**)$$
 If we call $d := (q_i - 1) / e$, which equals $b_i - 1$ except in the case quadratic crystalline case where $b_i = q_i$, $e=2$ and $d = (b_i - 1) / 2$, we have:
 $$  t_1 \equiv a \pmod{d}  \qquad \qquad (i) $$
 $$  t_2 \equiv a+1 \pmod{d} \qquad \qquad (ii) $$
 $$ 2a+1 \equiv k+d-1 \pmod{q_i - 1} \qquad \qquad (iii) $$
 Congruence (iii) just says that the determinant of $\bar{\rho}$ is, locally at $q_i$, equal to $\chi^{k+d-1}$: this follows from the fact that $\rho$ has determinant, locally at $q_i$,  equal to $\chi^{k-1} \cdot \omega_{q_i}^d$. \\
 Adding and subtracting equations (i) and (ii) and using (iii), we get (recall that $d$
divides $q_i-1$):
$$ t_1 + t_2 \equiv k-1 \pmod{d} $$
$$ t_2 - t_1 \equiv 1  \pmod{d}   $$
Now, inequalities (**) on $t_1$ and $t_2$ and the fact that $k \leq q_{i-1} + 1 < d$ show that these two congruences are
indeed equalities. Hence $t_1 = k/2 - 1$ and $t_2 = k/2$. Doubling (i), we
obtain $2a \equiv k-2 \pmod{2d}$ and then combining with (iii), we finally get
$ d \equiv 0 \pmod{2d}$ since $2d$ divides $q_i-1$ (recall that we know that $(q_i -1) / d = e$ is always
even). This is a contradiction, so we conclude that $\bar{\rho}_{q_i}$ satisfies Kisin's condition (T) in Theorem \ref{teo:superKisin}. 
\qed
\\

Remark: The good-dihedral prime $q$ is  the only prime that remains at the end in the level, then the bound $B$ such that $q$ is good-dihedral with respect to this bound (see section 2) must be chosen  larger than all primes $b_i$ and $q_i$ in the above construction, thus ensuring non-solvable residual image through the whole  process. To be more precise: given a residual representation in characteristic $p$ of odd conductor $N$ and weight $k$, before starting the ``iterated killing ramification" step one has to ``precompute" the primes $b_i$ and $q_i$ defined as in the previous proof, then fix a bound $B$ larger than all of them (thus in particular also larger than $p, k$ and the prime factors of $N$) and add to the level a prime $q$ that is good-dihedral with respect to this bound $B$. \\

         \section{Reduction to the level $3$ case}
         
         It follows from lemma \ref{teo:good-dihedral} and all the lemmas and the proposition in the previous section that we
have reduced the proof of the odd level case of Serre's conjecture to the case: characteristic $p$, weight $k$, level $q^2$, good dihedral at $q$, where $k$ and the characteristic $p$ are smaller that the bound $B$.\\


         \subsection{Weight reduction}
         
         \begin{lema}
         \label{teo:wr} Let $\bar{\rho}_p$ be a Galois representation of Serre's weight $k < B$, with $p<B$, and Serre's level $q^2$, and suppose that $\bar{\rho}_p$ is locally good-dihedral at $q$. Then, the modularity of $\bar{\rho}_p$ follows from the modularity of a representation $\bar{\rho}'_{p'}$ of level $3 \cdot q^2$ or $q^2$, locally good-dihedral at $q$, semistable at $3$ in the first case, of weight $2$.
         \end{lema}
         
 \bf{Proof}: \rm        By applying Khare's weight reduction, i.e., the method of weight reduction used in [Kh] to prove the level $1$ case of Serre's conjecture, we are reduced to the cases $k=2,4,6$, level $q^2$. The only difference with the case considered in [Kh] is that now the level is $q^2$, but since we are going to work in characteristics smaller than $q$ (moreover, smaller than $B$) we can perform weight reduction exactly as in the level $1$ case. Observe that all residual Serre's weights are going to be even because the ramification at $q$ has trivial character (i.e., the determinants of the representations are not ramified at $q$), and using the good-dihedral prime $q$ the considerations that are required in the level $1$ case to apply M.L.T. in residually reducible cases and other cases of small residual image are not going to be needed here. Finally, observe that since the result of ``existence of weight $2$ minimal lifts" that we stated in section 2 (taken from [KW2]) is stronger than the one used in [Kh], there is no need to treat the ordinary and the non-ordinary cases separately as in [Kh].\\
         
          Now we have to perform one more weight reduction in order to reduce the case $k=6$ to the cases $k=2$ or $4$. Recall that since we have the good-dihedral prime in the level residual images will be non-solvable.\\
         In this weight reduction modularity is going to be preserved due to the results in [Ki1] and [Ki3], and those in [SW2].  The trick that we will apply
 is the same used with the pair of primes $3$ and $2$ in [KW2] to reduce the weight $4$ and weight $3$ cases to the weight $2$ case (but here we do not
  need $2$-adic modularity lifting theorems since we do not work with $p=2$): We start with a residual representation of weight $6$, level $q^2$, which is good-dihedral at a large prime $q$. We switch to characteristic $5$, reduce mod $5$ (here we assume  that the residual Serre's weight is $6$: if not, it is equal to $2$ and thus we are done) and consider a weight $2$ lift, corresponding to an abelian variety of $\GL_2$-type with bad semistable reduction at $5$ and conductor of the corresponding system of $2$-dimensional Galois representations equal to $5 \cdot q^2$. Then we switch to characteristic $3$, reduce mod $3$ and here we observe that this mod $3$ representation, since it is either unramified or has unipotent ramification at $5$ (and in both cases we know that there is a lift with semistable ramification at $5$,
 so in the unramified case the well-known necessary condition for raising the level at $5$ is satisfied), and $3 \mid (5+1)$, admits a weight $2$ lift where the ramification at $5$ is no longer semistable but instead is given by a character of order $3$. This ``level raising" trick, as we have already remarked, is similar to the one used in [KW2], and the existence of such non-minimal weight $2$ lifts and the strongly compatible systems containing them is stated in [KW2], Theorem 5.1 (4) and proved in [KW3]. We obtain a lift of conductor $25 \cdot q^2$: what we have just constructed is a non-minimal lift (it is not minimal at $5$),
  having the same kind of ramification at $5$ and at the good-dihedral prime $q$. \\
  We consider the strongly compatible system containing this $3$-adic representation and we switch to characteristic $5$. Using strong compatibility and the description of ramification at $5$ (a character of order $3 \mid (5+1)$) it follows from results of D. Savitt (cf. [Sa]) that the residual mod $5$ representation will have (after suitable twist) Serre's weight equal to $2$ or $4$, but never $6$. This concludes the weight reduction.\\
  
   Now we take a crystalline minimal lift, switch to characteristic $3$ and reduce mod $3$.\\
       By considering a crystalline minimal lift (if $k=2$) or a minimal weight $2$ lift
         (if $k=4$) and reducing modulo some characteristic bigger than $3$ we see that we have reached a case of weight $2$ and level $q^2$ or $3q^2$, semistable at $3$.
        \qed 
         \\

                  Remark 1: We have used in the above proof  {\it minimal weight $2$ lifts for the case of Serre's weight $k= p+1$}.  \\
                  Remark 2: The modularity lifting result for a crystalline representation of weight $k$ (i.e., of Hodge-Tate weights $0$ and $k-1$) in characteristic $p = k-1$ is just a combination of the results of Skinner-Wiles in [SW2], for the ordinary case, with the main theorem of [Ki3], which works in the complementary case as follows from the description of crystalline non-ordinary representations given in [BLZ]. It does not require any condition on the residual image (it can even be reducible, a case where modularity follows from the results in [SW1]). Observe that in this situation the residual Serre's weight can be either $k'= p +1$  or $k'=2$ (cf. [BLZ]), in the above proof for simplicity we have always assumed that we are in the first case, which is the worst case, the second case leading (by a similar, yet simpler, argument) to the same conclusion. \\

         \subsection{Removing the good-dihedral prime}
         
         Starting from a residual representation as in the conclusion of the previous lemma, we take a minimal lift and insert it in a strongly compatible system. With this weight $2$ system of conductor $q^2$ or $3q^2$ we move to characteristic $t$ ($t$ as in the definition of good-dihedral prime given in section 2) and we consider the residual mod $t$ representation. This is the first time since the first step (when we added the good-dihedral prime $q$ to the level) that we have switched to a characteristic greater than the bound $B$.\\
          Moreover, since the character giving ramification at $q$ of our compatible system is a character of order $t^\alpha$, it is clear that we are losing some ramification at $q$ at this step, so a few remarks before going on:\\
         -Remark 3: This residual representation is no longer locally good-dihedral at the prime $q$. This means that this mod $t$ representation and any of the residual representations in the next steps may have solvable images or even be reducible, in which case we can not apply any longer the result of [Ki4]. We will have to use (see, for example, Khare's proof of the level $1$ case for a similar situation) other modularity lifting results. As already explained in the introduction (this is proved in [Di1], [DM] and [Kh]) in several cases this is known to work well (even if we have no information on the residual image):\\
         -  semistable weight $2$ lift\\
         -  potentially Barsotti-Tate lift which is Barsotti-Tate over the cyclotomic extension $\Q_p (\mu_p)$\\
         -  crystalline lift of weight $k$ with $k<p$ or $k=p+1$ (for $k=p+1$ either the lift is ordinary and the results of Skinner-Wiles in [SW1] and [SW2] apply or if not the result of Kisin in [Ki3] apply), assuming that $p \neq 2k-3$ or the residual representation is semistable at all primes different from $p$.\\
         In all the steps that follow we will always be in one of the above situations, thus ensuring that modularity propagates.\\
         -Remark 4: When reducing mod $t$, given the information on the ramification at $q$ of the system, there are two possibilities (this is also noticed in [KW2]): the residual representation is either unramified or semistable (i.e., unipotent ramification) at $q$. \\

         From now on, at each ``move" the following remark applies: if the residual representation has solvable image (or reducible) then it is modular (or reducible) and it can be seen that the given lift is modular (because of the above remark 3), therefore the proof of modularity concludes. So the only case relevant is the case where the residual image is non-solvable, which is the case we will consider.\\

         So, if the mod $t$ representation is unramified at $q$, it will have $k=2$ (the $t$-adic lift was Barsotti-Tate) and $N=3$ (semistable), a case of Serre's conjecture already solved (cf. [Di3], [KW1]). Thus, the case that remains is the case: $k=2$, $N= 3q$, semistable at both primes (the case $k=2$, $N=q$ has already been solved by Khare in [Kh]).\\
         We take a crystalline minimal lift (which corresponds to an abelian variety with semistable reduction at $3$ and $q$) and move to characteristic $q$ and reduce mod $q$. Since the residual weight will be either $2$ or $q+1$ we see that it only remains to solve the case $k = q+1$, $N=3$, semistable at $3$, i.e., the ``level $3$ case". \\
         Putting together lemma \ref{teo:good-dihedral}, the lemmas and the proposition in section 4, lemma \ref{teo:wr} and the above argument, we deduce the following:
         
         \begin{proposition}
         \label{teo:bajandohastaeltres} The truth of Serre's conjecture for any $\bar{\rho}_p$ ($p$ odd) of odd conductor follows from the truth of Serre's conjecture for all $\bar{\rho}'_{p'}$ of Serre's level $3$, semistable at $3$.
         \end{proposition}

    \section{Proof of the level $3$ case}

    As stated in the previous proposition, all that remains is to solve the level $3$ case (semistable at $3$). We apply Khare's weight reduction (as in [Kh]) and the proof is reduced to the cases: $k=2,4,6$, $N=3$ (semistable at $3$). This weight reduction is performed exactly as in [Kh], and in particular it follows as in [Kh] that in cases where residual representations are reducible or have small images a suitable M.L.T. still can be applied (see Remark 3 in the previous section, we are using the fact that the unique prime in the level is a semistable prime). The only slight difference with the weight reduction in the level $1$ case is the following: if, at some step of the weight reduction, one has to use as pivot the prime $3$ (to modify a character giving the ramification at $p$ for a prime $p$ such that $p-1$ is divisible by $3$, see [Kh]), then one proceed as in [Kh] moving to characteristic $3$ with a weight $2$ family, reducing mod $3$ (usually, the residual Serre's weight will be $4$, observe also that  M.L.T. work well in this semistable  weight $2$ situation, regardless of the size of the residual image: see Remark 3 in the previous section) and then taking a lift of this mod $3$ representation that is non-minimal at $p$ and is of weight $2$ (usually, this lift will be semistable at $3$): such lifts exist by [KW2], Theorem 5.1 (3). Of course, at some step (for example, at the one just described) the prime $3$ can disappear from the conductor, thus reducing the proof to the level $1$ case. \\

     Therefore, after applying Khare's weight reduction, since the case $k=2$ is known,  we only have to consider two cases.\\

   $\bullet$ Case $k=6$, $N=3$: Take a crystalline minimal lift and move to characteristic $5$, and reduce mod $5$. Since the case $(k,N) = (2,3)$ is known, we assume we are in the case $(k,N)=(6,3)$. Take a minimal weight $2$ lift, which is known to correspond, as follows from the potential modularity results of Taylor (cf. [Ta1] and [Ta2]), to a semistable abelian variety of $\GL_2$-type with good reduction outside $3$ and $5$. By recent results of Schoof (cf. [Sc2], the reader may also consult [Sc1] where a similar result is proved in the case of one single small prime of bad semistable reduction) such an abelian variety is known to be modular. This concludes the proof in this case.\\

    $\bullet$ Case  $k=4$, $N=3$: Take a crystalline minimal lift and move to characteristic $5$, and reduce mod $5$ (we assume that we are not losing ramification at $3$, if not, this is a level $1$ case which is already solved). If this residual representation is reducible or has solvable image, then modularity lifting theorems are enough to conclude the proof of this case of Serre's conjecture. On the other hand, if we assume that this residual representation has non-solvable image, we can take a minimal weight $2$ lift of it, whose inertial Weil-Deligne parameter at $5$ is $(\omega_5^2  \oplus 1,0)$. Such a lift is known to correspond (as follows again from the potential modularity results of Taylor) to an abelian variety $A$ of $\GL_2$-type with bad reduction at $3$ and $5$: semistable reduction at $3$ and potentially good reduction at $5$. Since ramification at $5$ is given by the quadratic character corresponding to $\Q(\sqrt{5})$ we know that $A$ acquires good reduction at $5$ over this real quadratic field, but this contradicts the other main result in [Sc2] which says that there are no semistable abelian varieties with good reduction outside $3$ over $\Q(\sqrt{5})$.\\ 
    
    This concludes the proof of the level $3$ case and, because of proposition
     \ref{teo:bajandohastaeltres}, of the odd level case of Serre's conjecture.\\

    Let us write the theorem we have proved, together with some well-known consequences (cf. [Se] for the proof of the second consequence, and [Ri1] for the proof of the first consequence).

\begin{teo}
\label{teo:SerreNK}
Serre's conjecture is true for any odd, two-dimensional, irreducible Galois representation whose Serre's level is odd.\\
Every abelian variety defined over $\Q$ of $\GL_2$-type having good reduction at $2$ is modular.\\
Every rigid Calabi-Yau threefold defined over $\Q$ having good reduction at $2$ is modular.
\end{teo}

\section{Final Remarks}

\subsection{On even levels}

a) There is one case of even level that can be proved by these methods: the case of level $6$ (semistable) and weight $2$. The method is the following: move to characteristic $3$, then the residual mod $3$ representation has conductor $2$ and weight $2$ or $4$. For these two cases, Serre's conjecture has been proved by Moon and Taguchi (cf. [MT], they proved reducibility, of course) in characteristic $3$, thus applying modularity lifting results (modularity of the semistable weight $2$ deformation) the proof is complete.\\
b) Assuming GRH, also the case of level $10$ (semistable) and weight $2$ is known: the crystalline minimal lift corresponds to a semistable abelian variety, which is modular by results of Calegari (assuming GRH, cf. [Cal]). Therefore, assuming GRH, we can also prove the following cases of Serre's conjecture: level $2p$ (semistable), $p$ any odd prime, weight $2$. The method is the following: move to characteristic $p$, then the proof is reduced to prove the case: level $2$, weight $k \geq 2$. Applying Khare's weight reduction, this can be solved for arbitrary weight assuming that some base cases are known: $k=2,4,6$, $N=2$. For $k=2,4$ this is known, thanks to the result of Moon and Taguchi. For $k=6$, $N=2$, we move to characteristic $5$ and the residual mod $5$ representation has a weight $2$ lift corresponding to a semistable abelian variety of conductor (dividing) $10$. Since assuming GRH such a variety is modular, we conclude the proof. As a corollary, it follows that any semistable abelian variety of $\GL_2$ type with bad reduction only at $2$ and an odd prime $p$ is modular, assuming GRH.

\begin{teo}
\label{teo:Serre22p} Assume the Generalized Riemann Hypothesis. Then
Serre's conjecture is true for any odd, two-dimensional, irreducible Galois representation of Serre's weight $2$ and semistable level $2p$
 for any odd prime $p$ (and the case $p=3$ holds unconditionally, i.e., independently of GRH).\\
Every semistable abelian variety defined over $\Q$ of $\GL_2$-type having bad reduction only at $2$ and another prime $p$ is modular.\\
\end{teo}

\subsection{Stronger versions of Kisin's results}

According to some experts' opinion, the modularity lifting result of [Ki4] could be improved, in particular it is expected that a proof of a stronger version without condition (T) should be given in the near future. It is an easy exercise to see that assuming that such a strong modularity lifting result holds our proof can be simplified, in particular the iterated killing ramification process in section 4. Also, the reduction to weights $k=2$ and $4$ at the beginning of section 5 can be obtained automatically by just switching to characteristic $3$. \\

\section{Bibliography}

[BLZ] Berger, L., Li, H., Zhu, H., {\it Construction of some families of $2$-dimensional crystalline representations}, Math. Annalen {\bf 329} (2004) 365-377
\newline
[Cal] Calegari, F., {\it Semistable Abelian Varieties over $\Q$},  Manuscripta Math. {\bf 113} (2004) 507-529
\newline
[Car] Caruso, X., {\it Repr\'{e}sentations semi-stables de torsion
 dans le cas $er < p - 1$}, J. Reine Angew. Math. {\bf 594} (2006) 35-92
\newline
[Co] Coleman, M. D., {\it On the equation $b_1 p - b_2 P_2 = b_3$}, J. Reine Angew. Math. {\bf 403} (1990) 1-66
\newline
[Dia] Diamond, F., {\it On deformation rings and Hecke rings},  Ann. of Math. {\bf 144} (1996) 137–166  
\newline
[DFG] Diamond, F., Flach, M., Guo, L., {\it The Tamagawa number conjecture of adjoint motives of modular forms}, Ann. Sci. Ec. Norm. Sup. {\bf 37} (2004) 663-727
\newline
[Di1] Dieulefait, L., {\it Modularity of Abelian Surfaces with Quaternionic
Multiplication}, Math. Res. Letters {\bf 10} n. 2-3 (2003)
\newline
[Di2] Dieulefait, L., {\it Existence of compatible families and new cases of the Fontaine-Mazur conjecture},
J. Reine Angew. Math. {\bf 577} (2004) 147-151
\newline
[Di3] Dieulefait, L., {\it The level $1$ weight $2$ case of Serre's conjecture}, Rev. Mat. Iberoamericana {\bf 23} (2007) 1115-1124
\newline
[Di4] Dieulefait, L., {\it Langlands base change for $\GL_2$}, preprint, (2009); available at: www.arxiv.org
\newline
[DM] Dieulefait, L., Manoharmayum, J., {\it Modularity of rigid Calabi-Yau threefolds over $\Q$}, in ``Calabi-Yau Varieties and Mirror Symmetry" (eds. Yui, N., Lewis, J.), Fields Inst. Comm. Series {\bf 38} (2003) 159-166
\newline
[Kh] Khare, C., {\it Serre's modularity conjecture: the level $1$ case}, Duke Math. J. {\bf 134} (2006)  557-589
\newline
[KW1] Khare, C., Wintenberger, J-P., {\it On Serre's reciprocity conjecture for $2$-dimensional
 mod $p$ representations of $\Gal(\bar{\Q}/\Q)$}, Annals of Math. {\bf 169} (2009) 229-253
 \newline
[KW2] Khare, C., Wintenberger, J-P., {\it Serre's modularity conjecture (I)}, Invent. Math. {\bf 178}   (2009) 485-504
\newline
[KW3] Khare, C., Wintenberger, J-P., {\it Serre's modularity conjecture (II)}, Invent. Math. {\bf 178}   (2009) 505-586
\newline
[Ki1] Kisin, M., {\it Moduli of finite flat group schemes, and modularity}, Annals of Math. {\bf 170} (2009) 1085-1180 
\newline
[Ki2] Kisin, M., {\it Modularity of $2$-dimensional Galois representations},
 Current Developments in Mathematics (2005) 191-230
\newline
[Ki3] Kisin, M., {\it Modularity of some geometric Galois representations},
 in ``L-functions and Galois representations" (Durham 2004), Cambridge U.P.  (2008) 438-470
\newline
[Ki4] Kisin, M., {\it The Fontaine-Mazur conjecture for $\GL_2$}, J. Amer. Math. Soc. {\bf 22} (2009) 641-690
\newline
[MT] Moon, H. ; Taguchi, Y.,  {\it Refinement of Tate's discriminant bound and non-existence theorems for mod $p$ Galois representations},  Documenta Math., Extra Volume: Kazuya Kato's Fiftieth Birthday (2003) 641-654
\newline
[Ri1] Ribet, K., {\it Abelian varieties over $\Q$ and modular forms}, Algebra and topology 1992 (Taej{\'o}n), 53-79, Korea Adv. Inst. Sci. Tech., Taej{\'o}n (1992)
\newline
[Ri2] Ribet, K., {\it Images of semistable Galois
representations}, Pacific J.  Math. {\bf 181} (1997)
\newline
[Sa] Savitt, D., {\it On a conjecture of Conrad, Diamond, and Taylor}, Duke Math. J. {\bf 128} (2005) 141-197
\newline
[Sc1] Schoof, R., {\it Abelian varieties over $\Q$ with bad reduction in one prime only}, Comp. Math. {\bf 141} (2005) 847-868
\newline
[Sc2] Schoof, R., {\it Semistable abelian varieties with good reduction outside $15$}, preprint (2009); available at: http://www.mat.uniroma2.it/$\sim$schoof/papers.html
\newline
[Se] Serre, J-P., {\it Sur les repr{\'e}sentations modulaires de degr{\'e}
$2$ de $\Gal(\bar{\mathbb{Q}} / \mathbb{Q})$}, Duke Math. J. {\bf 54}
(1987) 179-230
\newline
[SW1] Skinner, C., Wiles, A., {\it Residually reducible representations and modular forms}, Publ. Math. IHES {\bf 89} (1999) 5-126
\newline
[SW2] Skinner, C., Wiles, A., {\it Nearly ordinary deformations of irreducible
residual representations} Ann. Fac. Sci. Toulouse Math. (6) {\bf 10} (2001)      185-215
\newline
[Ta1] Taylor, R., {\it Remarks on a conjecture of Fontaine and Mazur},
J. Inst. Math. Jussieu {\bf 1} (2002)
\newline
[Ta2] Taylor, R., {\it On the meromorphic continuation of degree two
 L-functions}, Documenta Mathematica, Extra Volume: John Coates' Sixtieth Birthday (2006) 729-779
\newline
 [Xi] Xie, S. G., {\it On the equation $a p - b P_2 = m$},  Acta Math. Sinica {\bf 3} (1987) 54-57

\end{document}